\documentclass[3p,10pt,a4paper,twoside,fleqn,sort&compress]{filomat}
\usepackage{amssymb,amsmath,latexsym}
\usepackage[varg]{pxfonts}
\usepackage[toc,page]{appendix}
\usepackage{graphicx}
\newcommand{\FilVol}{33:7}
\newcommand{\FilYear}{2019}
\newcommand{\FilPages}{2135--2148}
\newcommand{\FilDOI}{https://doi.org/10.2298/FIL1907135C}
\newcommand{\FilReceived}{Received: 03 April 2017; Revised: 28 July 2017; Accepted: 10 January 2019}

\begin{document}

\title{A comparison between MMAE and SCEM for solving singularly perturbed linear boundary layer problems}

\author{S\"uleyman Cengizci}
\ead{suleyman.cengizci@antalya.edu.tr}
\address{Computer Programming, Antalya Bilim University, 07190, Antalya, Turkey}
\newcommand{\AuthorNames}{S\"uleyman Cengizci}

\newcommand{\FilMSC}{Primary 34B15, 34E15}
\newcommand{\FilKeywords}{Boundary layer, Singular perturbation, Successive complementary expansion method, Uniformly valid approximation}
\newcommand{\FilCommunicated}{Marko Petkovic}

\begin{abstract}
In this study, we propose an efficient method so-called \textit{Successive
Complementary Expansion Method (SCEM),} that is based on generalized
asymptotic expansions, for approximating to the solutions of singularly
perturbed two-point boundary value problems. In this easy-applicable method, in contrast to the well-known method \textit{the Method
of Matched Asymptotic Expansions (MMAE),} any matching process is not
required to obtain uniformly valid approximations. The key point: A
uniformly valid approximation is adopted first, and complementary functions
are obtained imposing the corresponding boundary conditions. An illustrative
and two numerical experiments are provided to show the implementation and
numerical properties of the present method. Furthermore, MMAE results are
also obtained in order to compare the numerical robustnesses of the methods.
\end{abstract}

\maketitle

\makeatletter
\renewcommand\@makefnmark%
{\mbox{\textsuperscript{\normalfont\@thefnmark)}}}
\makeatother

\section{Introduction}

Many phenomena in biology, chemistry, physics and engineering sciences are
modelled and formulated by boundary value problems associated with different
kinds of differential equations. In this manner, a model which is formulated
by an equation or a system containing positive small parameter(s) is
referred to as \emph{perturbed model} and a model which does not keep the
positive small parameter(s) is named as \emph{reduced} or \emph{unperturbed
model} \cite{1}. If the perturbed model contains the small parameter(s) as
coefficient(s) to the highest order derivative term(s), then the problem is
referred to as \emph{singularly perturbed problem}, otherwise called
as \emph{regularly perturbed problem}.

First studies on perturbation problems were conducted by L. Prandtl and J.
H. Poincar\'{e} in 1900's. The term \emph{"boundary layer"} first appeared
in Prandtl's paper \emph{"Motion of fluids with very little viscosity"} in
1904. In this work, the small parameter was the inverse Reynold number and
the equations were based on the classical Navier-Stokes equations of fluid
mechanics \cite{2}. Since 1900's, various methods have been constructed and
employed. Kadalbajoo and Reddy \cite{3} made a survey of various asymptotic
and numerical methods that were developed between 1908 -1986 for the
determination of approximate solutions of singular perturbation problems,
then Kadalbajoo and Patidar \cite{4} extended the work done by Kadalbajoo and
Reddy and they surveyed the studies done by various researchers in the field of
singular perturbation problems between 1984-2000 considering only one
dimensional problems and made discussion on linear, nonlinear,
semilinear and quasi-linear problems. Subsequently, Kadalbajoo and Gupta \cite{5}
presented a great survey on computational methods for different types of
singular perturbation problems solved by different researchers between
2000-2009. Kumar studied general methods for solving singularly perturbed
problems arising in engineering in his work \cite{6}. Later in \cite{7}, Roos
made a survey in singularly perturbed partial differential equations
covering the years 2008-2012. Besides these great papers and surveys, there
are also great reference books on applications and theory of singular perturbation
problems. Some of them may be given as by Van Dyke\cite{8}, Holmes \cite{9},
Johnson \cite{10}, Lagerstrom \cite{11}, O'Malley \cite{12}, Eckhaus \cite{13},
Verhulst \cite{14}, Paulsen \cite{15}, Nayfeh \cite{16}, Murdock \cite{17},
Murray \cite{18}, and by Hinch \cite{19}.

In this paper, we employ the MMAE to approximate to the solutions of
singular perturbation problems first, and later present an efficient
asymptotic method which was introduced by Jean Cousteix and Jacques Mauss in 
\cite{20,22} as an alternative method to the MMAE: Successive Complementary
Expansion Method (SCEM). The main principle of SCEM is built on the purpose
of obtaining uniformly valid approximation to the singular perturbation
problems without any matching procedure. The method is, in general,
applicable to singular perturbation problems that we can approximate by
MMAE. To this end, we can also compare the SCEM results with previously
obtained ones by MMAE. The first step in SCEM starts with looking for an
approximation for outer region. The approximation is generally in good
quality for the outer region, but not for the inner region. The main idea
is, using boundary conditions, to add a correction term that complements the
approximation. The procedure can be iterated using new corrections for new
terms to improve the accuracy of approximation. The most important advantage
of SCEM is its ability of giving uniformly valid approximation without any
matching procedure. Boundary conditions are enough to implement the method.
Moreover, the boundary conditions are satisfied exactly, but not
asymptotically.

\section{Description of the method}
In this work, we deal with singularly perturbed second order linear
two-point boundary value problems in the following form

\begin{equation*}
\varepsilon y^{\prime \prime }(x)+p(x)y^{\prime }(x)+q(x)y(x)=r(x),\text{ }%
a\leq x\leq b;\text{ }a,b\in 
\mathbb{R}%
\end{equation*}%
with the boundary conditions $y(a)=\alpha \in 
\mathbb{R}
$ and $y(b)=\beta \in 
\mathbb{R}
,$ where $0<\varepsilon \ll 1$, $p(x),q(x)$ and $r(x)$ are sufficiently\
smooth functions. As $\varepsilon \rightarrow 0^{+},$ the order of the
differential equation is reduced and the equation that we call \textit{%
reduced equation} 
\begin{equation*}
p(x)y^{\prime }(x)+q(x)y(x)=r(x),\text{ }a\leq x\leq b
\end{equation*}%
is formed. One can observe that there are two boundary conditions in the
original problem but only one of them can be imposed to the reduced
equation. Moreover, as $\varepsilon $ tends to $0,$ because of the reduction
of the order, rapid changes occur in the solution. The region in which these
rapid changes occur is named as \emph{inner layer} or \textit{boundary layer}, and the layer
in which the solution exhibits mildly changes is named as \emph{outer layer}%
. The sign of the coefficient function $p(x)$ determines the type of the
layer(s). Over the interval $\left[ a,b\right] ,$ if $p(x)>0$ for all $x$,
then a boundary layer occurs at the left-end of the interval, if $p(x)<0$
for all $x,$ then a boundary layer occurs at the right-end of the interval
and if $p(x)$ changes sign in $(a,b),$ then interior layer(s) occurs at the
zero(s) of $p(x).$

\subsection{Asymptotic Approximations}

Given two functions $\phi (x,\varepsilon )$ and $\phi _{a}(x,\varepsilon ),$
defined in a domain $\Omega ,$ are referred to as \emph{asymptotically
identical to the order }$\delta \left( \varepsilon \right) $ if their
difference is asymptotically smaller than $\delta \left( \varepsilon \right) $, 

\begin{equation}
\phi (x,\varepsilon )-\phi _{a}(x,\varepsilon )=o(\delta (\varepsilon )), 
\tag{2.1.1}
\end{equation}%
where $\delta \left( \varepsilon \right) $ is an order function and $\varepsilon $ is a positive small parameter arising from the physical
problem under consideration. The function $\phi _{a}(x,\varepsilon )$ is
named as an \emph{asymptotic approximation} of the function $\phi
(x,\varepsilon ).$

Asymptotic approximations, in general form, are defined by

\begin{equation}
\phi _{a}(x,\varepsilon )=\sum_{i=0}^{n}\delta _{i}(\varepsilon )\varphi
_{i}(x,\varepsilon ),  \tag{2.1.2}
\end{equation}%
where the asymptotic sequence of order functions $\delta _{i}(\varepsilon )$
satisfies the condition $\delta _{i+1}(\varepsilon )=o(\delta
_{i}(\varepsilon ))$, as $\varepsilon \rightarrow 0$. In these conditions,
the approximation (2.1.2) is named as \emph{generalized asymptotic expansion}%
. If the expansion\ (2.1.2) is written in the form of
\begin{equation}
\phi _{a}(x,\varepsilon )=E_{0}\phi =\sum_{i=0}^{n}\delta
_{i}^{(0)}(\varepsilon )\varphi _{i}^{(0)}(x),  \tag{2.1.3}
\end{equation}%
then it is named as \emph{regular asymptotic expansion}. The special
operator $E_{0}$ is called\textit{\ }\emph{outer\ expansion operator} at a
given order $\delta (\varepsilon )$, thus $\phi -E_{0}\phi =o(\delta
(\varepsilon ))$. For more detailed information on asymptotic
approximations, we refer the interested reader to \cite{9,12,13,18,20}.

\subsection{The MMAE for SCEM}

Interesting cases occur when the function is not regular in $\Omega ,$ so
one of the approximations (2.1.2) and (2.1.3) is valid only in a restricted
region $\Omega _{0}$ $\subset $ $\Omega $, that is called the \emph{outer
region}. Here, in the simplest case, we introduce
an inner domain which can be formally denoted as $\Omega _{1}=\Omega -\Omega
_{0}$ and located near the origin. In general, the boundary layer variable
is expressed as $\overline{x}=\frac{x-x_{0}}{\xi (\varepsilon )}$ where $x_{0%
\text{ }}$is the point at which the rapid changes begin to occur and $\xi
(\varepsilon )$ is the order of thickness of this boundary layer. If a
regular expansion can be constructed in $\Omega _{1}$, we can write

\begin{equation}
\phi _{a}(x,\varepsilon )=E_{1}\phi =\sum_{i=0}^{n}\delta
_{i}^{(1)}(\varepsilon )\varphi _{i}^{(1)}(\overline{x}),  \tag{2.2.1}
\end{equation}%
where the \emph{inner expansion operator} $E_{1}$ is defined in $\Omega _{1}$
of the same order $\delta (\varepsilon )$ just like the outer expansion
operator $E_{0};$ thus, $\phi -E_{1}\phi =o(\delta (\varepsilon ))$. As a
result,

\begin{equation}
\phi _{a}=E_{0}\phi +E_{1}\phi -E_{1}E_{0}\phi  \tag{2.2.2}
\end{equation}%
is clearly uniformly valid approximation\cite{11,13,21}.

In MMAE, two distinct approximations are found for two distinct (outer and
inner) regions and then to obtain a uniformly valid approximation over the
whole domain, these approximations are matched using limit process. Despite all the valuable works devoted to MMAE, it is not possible to
formulate a general mathematical theory of the method. Therefore, we will
study it on an illustrative example. 

Let us consider the following second order singularly perturbed problem that has been studied employing the MMAE in \cite{9}

\begin{equation}
\varepsilon y^{\shortmid \shortmid }+2y^{\shortmid }+2y=0,\text{ }y(0)=0,%
\text{ }y(1)=1,\text{ }0<x<1.  \tag{2.2.3}
\end{equation}%
This problem exhibits rapid changes near the point $x=0$ as $\varepsilon
\rightarrow 0^{+}$, this region is named as \textit{boundary layer}\ or 
\textit{inner layer}. But, over the other region which is called \textit{%
outer region},\ the solution does not show an unusual behavior. This region
is the region which is far from the point $x=0$. We will adopt an
approximation for the outer solution as

\begin{equation}
y(x)\approx y_{0}(x)+\varepsilon y_{1}(x)+\varepsilon ^{2}y_{2}(x)+...\text{
.}  \tag{2.2.4}
\end{equation}%
If approximation (2.2.4) is substituted into (2.2.3), we reach
\begin{equation}
\varepsilon \left( y_{0}^{\shortmid \shortmid }(x)+\varepsilon
y_{1}^{\shortmid \shortmid }(x)+...\right) +2\left( y_{0}^{\shortmid
}(x)+\varepsilon y_{1}^{\shortmid }(x)+...\right) +2\left(
y_{0}(x)+\varepsilon y_{1}(x)+...\right) =0,  \tag{2.2.5}
\end{equation}%
and letting $\varepsilon =0$, the reduced equation is obtained as%
\begin{equation}
y_{0}^{\shortmid }+y_{0}=0.  \tag{2.2.6}
\end{equation}%
Solution of (2.2.6) is obviously $y_{0}(x)=Ae^{-x},A\in 
\mathbb{R}
$ and if the outer boundary condition is imposed (for $x=1$), the outer
solution is obtained as

\begin{equation}
y_{0}(x)=e^{1-x}.  \tag{2.2.7}
\end{equation}%
In order to obtain inner (boundary) layer approximation, we introduce a new
variable

\begin{equation}
\overline{x}=\frac{x}{\varepsilon }\text{.}  \tag{2.2.8}
\end{equation}%
Thanks to the new variable, called \textit{boundary layer (stretching)
variable},\ we get the chance to stretch the thin layer as $\varepsilon
\rightarrow 0^{+}.$ We denote the inner solution that depends on$\ \overline{%
x}\ $and valid for near the point $x=0$ by $Y(\overline{x}).$ Using the
chain rule
\begin{equation}
\frac{d}{dx}=\frac{d\overline{x}}{dx}\frac{d}{d\overline{x}}=\frac{1}{%
\varepsilon }\frac{1}{d\overline{x}}  \tag{2.2.9}
\end{equation}%
is obtained and applying the transformation (2.2.8) to the original problem (2.2.3)
\begin{equation}
\varepsilon ^{-1}\frac{d^{2}Y}{d\overline{x}^{2}}+2\varepsilon ^{-1}\frac{dY%
}{d\overline{x}}+2Y=0  \tag{2.2.10}
\end{equation}%
is found. If the both sides of (2.2.10) are multiplied by $\varepsilon ,$\ we
reach to

\begin{equation}
\frac{d^{2}Y}{d\overline{x}^{2}}+2\frac{dY}{d\overline{x}}+2\varepsilon Y=0.
\tag{2.2.11}
\end{equation}%
Equation (2.2.11) is a regularly perturbed linear ordinary differential
equation and we propose an approximation in the form of

\begin{equation}
Y(\overline{x})\approx Y_{0}(\overline{x})+\varepsilon Y_{1}(\overline{x}%
)+\varepsilon ^{2}Y_{2}(\overline{x})+...\text{ .}  \tag{2.2.12}
\end{equation}%
Since we are interested only in the first term of the approximation
(2.2.12), we can make our calculations for $\varepsilon =0$. To this end,
the equation

\begin{equation}
\frac{d^{2}Y_{0}}{d\overline{x}^{2}}+2\frac{dY_{0}}{d\overline{x}}=0 
\tag{2.2.13}
\end{equation}%
is obtained. Considering the inner boundary condition at $x=0$ $%
(x=0\Longrightarrow \overline{x}=0),$ the general solution to the equation
(2.2.13) is found as

\begin{equation}
Y_{0}(\overline{x})=B\left(1-e^{-2\overline{x}}\right),B\in 
\mathbb{R}
.  \tag{2.2.14}
\end{equation}%
It is obvious that we are not able to determine one of the unknown constants, $B$.
To determine $B$, we will use the matching procedure of MMAE. So far we have
found two approximations (2.2.7) and (2.2.14) that are valid for distinct
regions, but we know that these two approximations actually belong to the
same approximation. Using this idea, the matching procedure may be given as
follows \cite{9,11}

\begin{equation}
\underset{x\rightarrow 0^{+}}{\lim }y_{0}(x)=\underset{x\rightarrow 0^{+}}{%
\lim }e^{1-x}=\underset{\overline{x}\rightarrow \infty }{\lim }B\left(1-e^{-2%
\overline{x}}\right)=\underset{\overline{x}\rightarrow \infty }{\lim }Y_{0}(%
\overline{x}).  \tag{2.2.15}
\end{equation}%
Thus, we reach

\begin{equation}
Y_{0}(\overline{x})=\left(e-e^{1-2\overline{x}}\right).  \tag{2.2.16}
\end{equation}%
Finally, we desire to obtain a composite approximation. To do this, we
simply add the inner and outer approximations and subtract the common limit,
that is, using the following procedure

\begin{align}
y & \approx y_{0}(x)+Y_{0}(\overline{x})-Y_{0}(\infty) \nonumber  \\
(\text{or equivalently }y & \approx y_{0}(x)+Y_{0}(\overline{x})-y_{0}(0^{+})) \tag{2.2.17}
\end{align}

we reach the following composite MMAE approximation%
\begin{equation}
y_{composite}\approx e^{1-x}-e^{1-\frac{2x}{\varepsilon }}.  \tag{2.2.18}
\end{equation}

\subsection{\textbf{Successive Complementary Expansion Method}}

The uniformly valid SCEM approximation is in the regular form given as

\begin{equation}
y_{n}^{scem}(x,\overline{x},\varepsilon )=\overset{n}{\underset{i=0}{\sum }}%
\delta _{i}(\varepsilon )\left[ y_{i}(x)+\Psi _{i}(\overline{x})\right] 
\tag{2.3.1}
\end{equation}%
where$\ \left\{ \delta _{i}\right\} \ $is an asymptotic sequence and
functions $\Psi _{i}(\overline{x})$ are the complementary functions that
depend on $\overline{x}$. Functions $y_{i}(x)$ are the outer approximation
functions that have been found by MMAE before, and they only depend on $x$,
not also on $\varepsilon $. If the functions $y_{i}(x)$ and $\Psi _{i}(%
\overline{x})$ also depend on $\varepsilon $, the uniformly valid SCEM
approximation, that is named as \emph{generalized SCEM approximation}\textit{%
, }is given in the following form \cite{22}-\cite{24}

\begin{equation}
y_{ng}^{scem}(x,\overline{x},\varepsilon )=\overset{n}{\underset{i=0}{\sum }}%
\delta _{i}(\varepsilon )\left[ y_{i}(x,\varepsilon )+\Psi _{i}(\overline{x}%
,\varepsilon )\right] .  \tag{2.3.2}
\end{equation}%
Let us consider the problem (2.2.3) again and propose an approximation for $%
n=0$, that is, we look for a SCEM approximation in the form of

\begin{equation}
y_{0}^{scem}(x,\overline{x},\varepsilon )=y_{0}(x,\varepsilon )+\Psi _{0}(%
\overline{x},\varepsilon )  \tag{2.3.3}
\end{equation}%
and we know that from the equation (2.2.7),\ $y_{0}(x)=e^{1-x}.$ Thus,
approximation (2.3.3) turns into

\begin{equation}
y_{0}^{scem}(x,\overline{x},\varepsilon )=e^{1-x}+\Psi _{0}(\overline{x}%
,\varepsilon ).  \tag{2.3.4}
\end{equation}%
Once the approximation (2.3.4) is substituted into the original problem
(2.2.3)

\begin{equation}
\varepsilon \frac{d^{2}}{dx^{2}}y_{0}^{scem}(x,\overline{x},\varepsilon )+2\frac{%
d}{dx}y_{0}^{scem}(x,\overline{x},\varepsilon )+2y_{0}^{scem}(x,\overline{x}%
,\varepsilon )=0  \tag{2.3.5}
\end{equation}%
is obtained. It follows that

\begin{align}
\varepsilon \frac{d^{2}}{dx^{2}}\left( e^{1-x}+\Psi _{0}(\overline{x}%
,\varepsilon )\right) +2\frac{d}{dx}\left( e^{1-x}+\Psi _{0}(\overline{x}%
,\varepsilon )\right) +2\left( e^{1-x}+\Psi _{0}(\overline{x},\varepsilon
)\right) & =0  \notag \\
\varepsilon e^{1-\varepsilon \overline{x}}+\varepsilon \frac{1}{\varepsilon
^{2}}\Psi _{0}^{\shortparallel }(\overline{x},\varepsilon
)-2e^{1-\varepsilon \overline{x}}+\frac{2}{\varepsilon }\Psi _{0}^{\shortmid
}(\overline{x},\varepsilon )+2e^{1-\varepsilon \overline{x}}+2\Psi _{0}(\overline{x},\varepsilon )& =0 
\notag \\
\frac{1}{\varepsilon }\Psi _{0}^{\shortparallel }(\overline{x},\varepsilon )+%
\frac{2}{\varepsilon }\Psi _{0}^{\shortmid }(\overline{x},\varepsilon
)+2\Psi _{0}(\overline{x},\varepsilon )& =-\varepsilon e^{1-\varepsilon 
\overline{x}}  \notag \\
\Psi _{0}^{\shortparallel }(\overline{x},\varepsilon )+2\Psi _{0}^{\shortmid
}(\overline{x},\varepsilon )+2\varepsilon \Psi _{0}(\overline{x},\varepsilon
)& =-\varepsilon ^{2}e^{1-\varepsilon \overline{x}}.  \tag{2.3.6}
\end{align}%
The resulting equation that is given in the last line of (2.3.6) is a
regularly perturbed linear non-homogeneous second order ordinary
differential equation. In order to obtain first term of SCEM approximation, letting $\varepsilon=0$, the first complementary SCEM approximation $\Psi _{0}(%
\overline{x},\varepsilon )$ is obtained as the solution of the following two-point boundary value problem 
\begin{equation}
\Psi _{0}^{\shortparallel }(\overline{x},\varepsilon )+2\Psi _{0}^{\shortmid
}(\overline{x},\varepsilon )=0,\ \Psi_{0} (0,0,\varepsilon )=-e,\ \Psi _{0} \Big(1,\frac{1}{\varepsilon },\varepsilon \Big)=0. \notag
\end{equation}
The first SCEM approximation is found as
\begin{equation*}
y_{0}^{scem}(x,\overline{x},\varepsilon )=e^{1-x}+\frac{e\left( e^{{-2\overline{x}}}-1\right) }{e^{\frac{-2}{\varepsilon }}-1}-e.
\end{equation*}
On the other hand, Problem (2.2.3) is an exactly solvable
problem and its exact solution is given as

\begin{equation}
y_{exact}(x)=\frac{\left( e^{\frac{x}{\varepsilon }(-1+\sqrt{1-2\varepsilon }%
)}-e^{\frac{x}{\varepsilon }(-1-\sqrt{1-2\varepsilon })}\right) }{e^{\frac{1%
}{\varepsilon }(-1+\sqrt{1-2\varepsilon })}-e^{\frac{1}{\varepsilon }(-1-%
\sqrt{1-2\varepsilon })}}.  \tag{2.3.7}
\end{equation}%

\begin{table}[ht]
\begin{center}
\begin{tabular}{|c|c|c|}
\hline
$\varepsilon $ & $L^{2}$ error in MMAE & $L^{2}$ error in SCEM \\ 
\hline
$0.0001$ & $0.000624687980610$ & $0.000624687980610$ \\ 
$0.0005$ & $0.003124942983068$ & $0.003124942983068$ \\ 
$0.0010$ & $0.006253648308213$ & $0.006253648308213$ \\ 
$0.0050$ & $0.031287231692987$ & $0.031287231692987$ \\ 
$0.0100$ & $0.061951928162705$ & $0.061951928162705$ \\ 
$0.0500$ & $0.283794475853395$ & $0.283794475853395$ \\ 
$0.1000$ & $0.498739448296190$ & $0.498739403531008$ \\ 
$0.3000$ & $0.609196858399138$ & $0.582659431139973$ \\ 
$0.4000$ & $0.540710242114810$ & $0.416122748110290$ \\ 
$0.6000$ & $0.893536533505815$ & $0.302123969698791$ \\ 
$0.8000$ & $1.754202335976711$ & $0.240771410186466$ \\ 
$1.0000$ & $2.790418350467303$ & $0.212524676097187$ \\
\hline
\end{tabular}
\end{center}
\caption{$L^{2}$-norm errors in MMAE and SCEM approximations for Illustrative Example}
\label{table:1}
\end{table}

\begin{figure}
    \centering
    \includegraphics[scale=0.4]{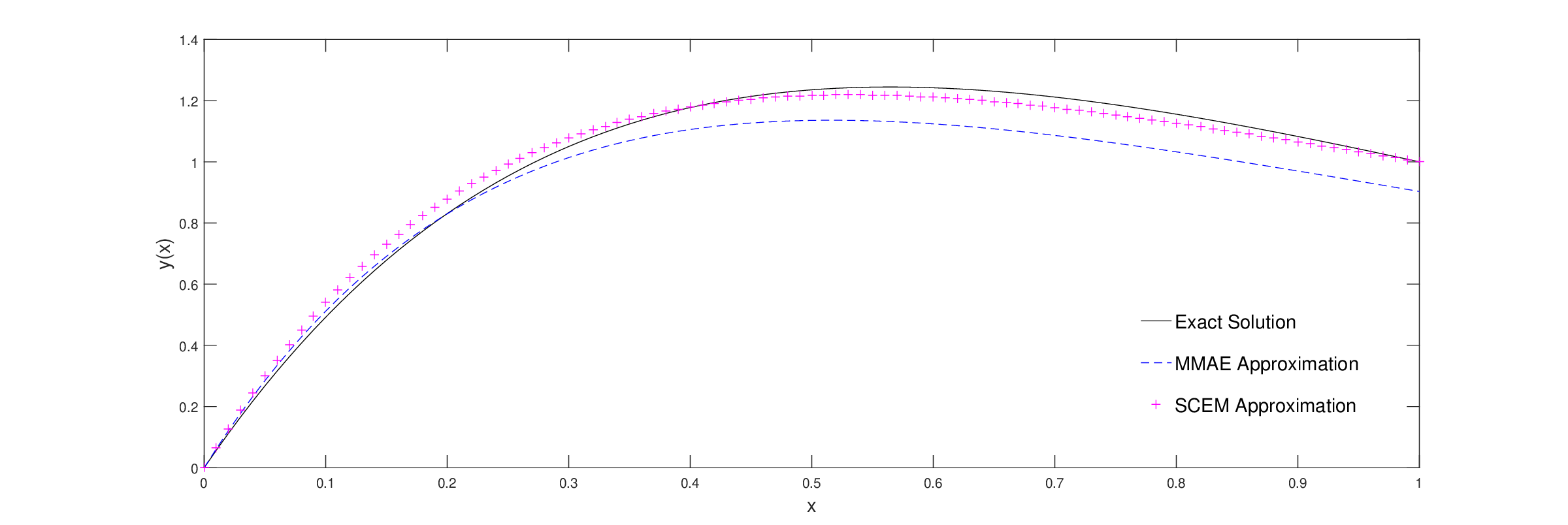}
    \caption{Comparison of results for illustrative example, $\protect\varepsilon =0.6$}
    \label{fig:my_label1}
\end{figure}

\begin{figure}
    \centering
    \includegraphics[scale=0.4]{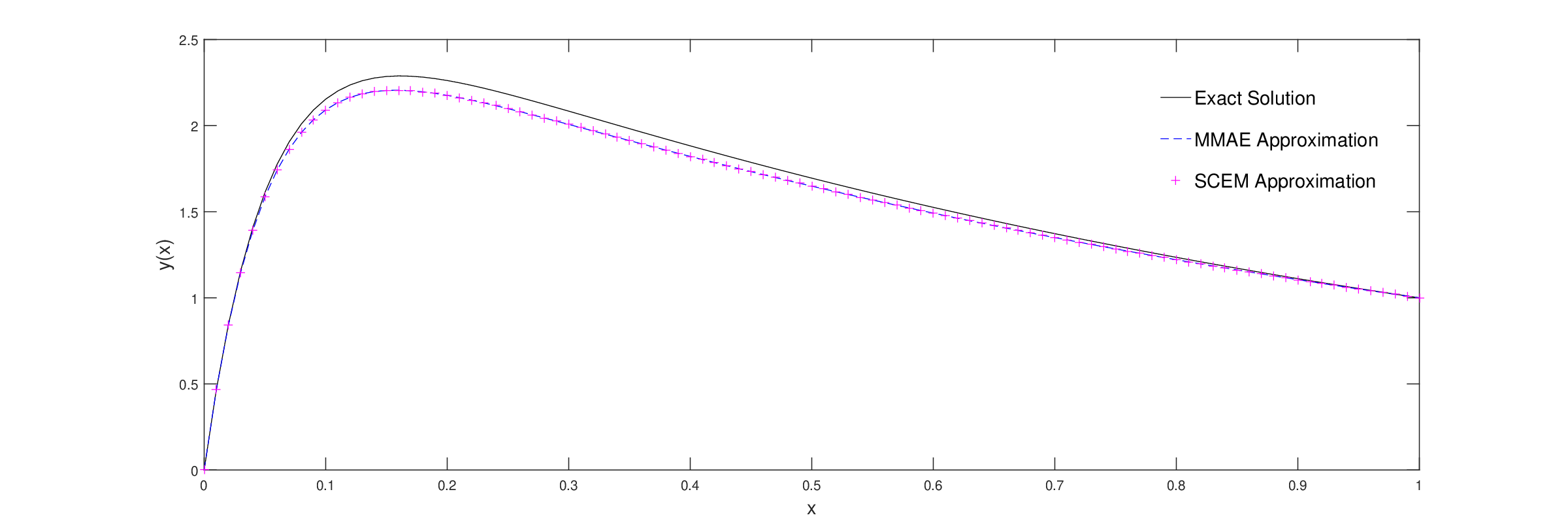}
    \caption{Comparison of results for illustrative example, $\protect\varepsilon =0.1$}
    \label{fig:my_label2}
\end{figure}

\begin{figure}
    \centering
    \includegraphics[scale=0.4]{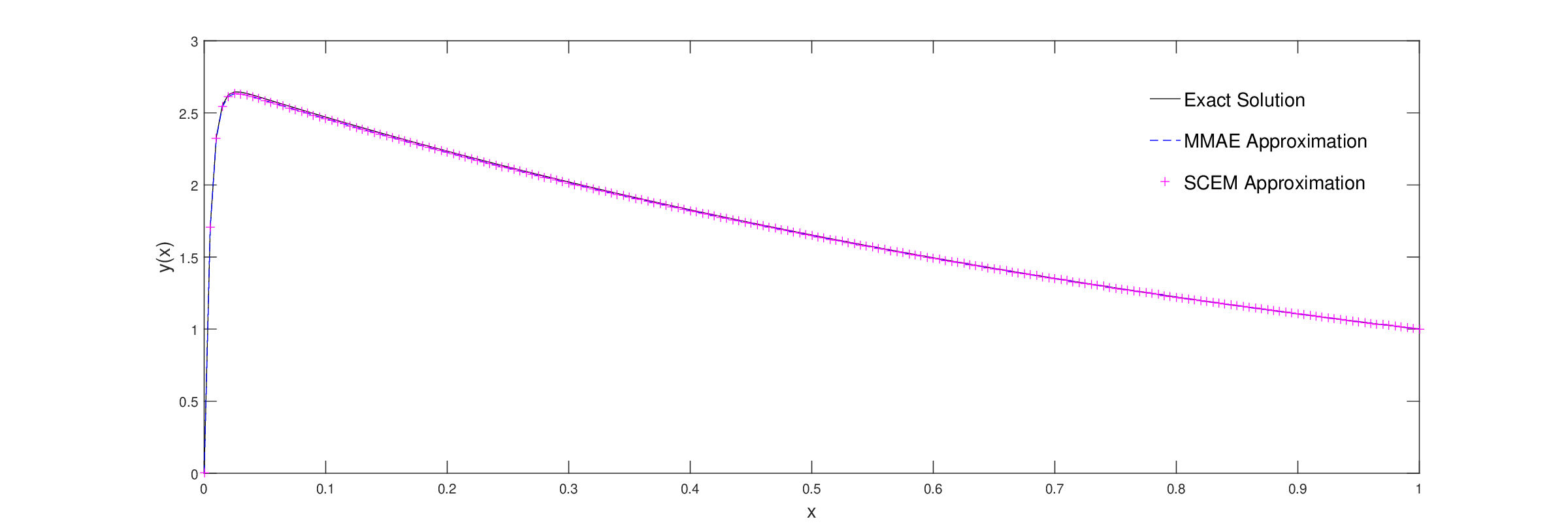}
    \caption{Comparison of results for illustrative example, $\protect\varepsilon =0.01$}
    \label{fig:my_label3}
\end{figure}
In Fig.\ref{fig:my_label1} - Fig.\ref{fig:my_label3}, exact solution, MMAE and SCEM approximations are given for the illustrative example. Here, the SCEM approximation represents the expression given by $y_{0}^{scem}(x,\overline{x},\varepsilon )=e^{1-x}+\frac{e\left( e^{\frac{-2x%
}{\varepsilon }}-1\right)}{e^{\frac{-2}{\varepsilon }}-1}-e$, that is obtained by balancing the equation (2.3.6) with respect to zeroth power (dominant order) of $\varepsilon$. While SCEM gives more accurate approximation for $\varepsilon=0.6$, as $\varepsilon$ gets smaller and smaller, SCEM and MMAE approximations overlap and approach to the exact solution. On the other hand, it can be observed especially in Fig.\ref{fig:my_label1} that MMAE approximations satisfy the boundary conditions asymptotically. Moreover, Table \ref{table:1} shows that the SCEM approximations are in quite good agreement with exact solution even increasing $\varepsilon$ values, and $L^{2}$-norm errors for these two approximations are exactly same up to the values around $\varepsilon=0.05$.   

\section{\textbf{Numerical Experiments}}

In this section, two numerical examples are given. All the figures and numerical results are generated in Matlab 2016b environment.

\textbf{Example 1:} Consider the non-homogeneous singularly perturbed
problem from fluid dynamics for fluid of small viscosity \cite{25}

\begin{equation}
\varepsilon y^{\shortmid \shortmid }+y^{\shortmid }=1+2x,\text{ }0\leq x<1,%
\text{ }y(0)=0,\text{ }y(1)=1.  \tag{3.1}
\end{equation}%
This problem has rapid changes near the point $x=0$. The exact solution to
(3.1) is given as

\begin{equation}
y_{exact}(x)=x^{2}+(1-2\varepsilon )x+\frac{2\varepsilon -1}{1-e^{\frac{-1}{%
\varepsilon }}}\left( 1-e^{\frac{-x}{\varepsilon }}\right)   \tag{3.2}
\end{equation}%
and one-term SCEM and MMAE approximations are obtained as follows.

The first term of the outer approximation, letting $\varepsilon =0$ , is obtained as the solution of initial value problem
\begin{align}
y_{0}^{\shortmid }=1+2x,\text{ }y_{0}(1)=1 \tag{3.3}
\end{align}
and it is clear that the solution is
\begin{equation}
    y_{0}(x)=x^{2}+x-1. \tag{3.4}
\end{equation}
Since the boundary layer occurs near the left-end of the interval, the
stretching variable should be in the form of $\overline{x}=\frac{x}{%
\varepsilon }$. The first term of the inner approximation, letting $\varepsilon =0$, is obtained as the solution of 
\begin{equation}
Y_{0}^{^{\shortparallel }}+Y_{0}^{^{\shortmid }}=0,\ Y_{0}(\overline{x}=0)=0\ \tag{3.5}
\end{equation}
and it is clear that the solution is
\begin{equation}
Y_{0}(\overline{x})=c_{1}(e^{-\overline{x}}-1),\ c_{1} \in \mathbb{R}. \tag{3.6}
\end{equation}
Now that the first terms of outer and inner approximations are obtained, these approximations should be matched by the definition of MMAE.
\begin{equation}
\underset{x\longrightarrow 0^{+}}{\lim }{}y_{0}{(x)}=\underset{x\longrightarrow 0^{+}}{\lim }(x^{2}+x-1)=
\underset{\overline{x}\rightarrow \infty }{%
\lim }c_{1}(e^{-\overline{x}}-1)=\underset{\overline{x}\longrightarrow \infty }{\lim }
Y_{0}(\overline{x}) \tag{3.7}
\end{equation}
Thus, it is determined that
\begin{equation}
Y_{0}(\overline{x})=e^{\overline{x}}-1. \tag{3.8}
\end{equation}
Finally, the composite MMAE approximation to the problem (3.1) is found as
\begin{equation}
y_{composite}(x,\varepsilon ) \approx x^{2}+x-1+e^{-\frac{x}{\varepsilon }}. \tag{3.9}
\end{equation}
Now, the one-term SCEM approximation can be determined substituting the approximation 
\begin{equation}
y_{0}^{scem}(x,\overline{x},\varepsilon )=y_{0}(x,\varepsilon )+\Psi _{0}(%
\overline{x},\varepsilon )=x^{2}+x-1+\Psi _{0}(%
\overline{x},\varepsilon ) \tag{3.10}
\end{equation}
into the original problem (3.1):

\begin{align}
\varepsilon \frac{d^2}{dx^2} \Big (x^{2}+x-1+\Psi_{0}(\overline{x},\varepsilon )\Big )+ \frac{d}{dx}\Big(x^{2}+x-1+\Psi_{0}(\overline{x},\varepsilon )\Big)& \notag =1+2\varepsilon\overline{x} \\ \notag
\varepsilon \Big (2+\frac{\Psi_{0}^{\shortparallel}(\overline{x},\varepsilon )}{\varepsilon^{2}}\Big )+ \Big( 2x+1+\frac{\Psi_{0} ^{\shortmid}(\overline{x},\varepsilon )}{\varepsilon}\Big)& =1+2\varepsilon\overline{x}\\ \tag{3.11}
2\varepsilon ^{2}+\Psi_{0}^{\shortparallel}(\overline{x},\varepsilon )+\varepsilon(2x+1) +\Psi_{0}^{\shortmid}(\overline{x},\varepsilon ) & =\varepsilon \big ( 1+2\varepsilon\overline{x} \big).
\end{align}
In order to obtain first term of SCEM approximation, letting $\varepsilon=0$, the first complementary SCEM approximation $\Psi _{0}(%
\overline{x},\varepsilon )$ is obtained as the solution of the following two-point boundary value problem 
\begin{equation}
\Psi_{0}^{\shortparallel}(\overline{x},\varepsilon )+\Psi_{0}^{\shortmid}(\overline{x},\varepsilon )=0,\ \Psi_{0} (0,0,\varepsilon )=1,\ \Psi_{0} \Big(1,\frac{1}{\varepsilon },\varepsilon \Big)=0. \tag{3.12}
\end{equation}
Finally, the first SCEM and MMAE approximation are found as follows
\begin{equation}
y_{0}^{scem}(x,\overline{x},\varepsilon )=x^{2}+x-1+\frac{e^{-\frac{1}{%
\varepsilon}}-e^{-\frac{x}{%
\varepsilon}}}{e^{-\frac{1}{\varepsilon }}-1}.  \tag{3.13}
\end{equation}
\begin{equation}
y_{composite}(x,\varepsilon ) \approx x^{2}+x-1+e^{-\frac{x}{\varepsilon }}. 
\tag{3.14}
\end{equation}

\begin{table}
\begin{center}
\begin{tabular}{|c|c|c|}
\hline
$\varepsilon $ & $L^{2}$ error in MMAE & $L^{2}$ error in SCEM \\ 
\hline
$0.0001$ & $0.001146036648629$ & $0.001146036648629$ \\ 
$0.0005$ & $0.005730183242787$ & $0.005730183242787$ \\ 
$0.0010$ & $0.011460350798498$ & $0.011460350798498$ \\ 
$0.0050$ & $0.057047561195172$ & $0.057047561195172$ \\ 
$0.0100$ & $0.112864481039688$ & $0.112864481039688$ \\ 
$0.0500$ & $0.513159514418249$ & $0.513159531588828$ \\ 
$0.1000$ & $0.901557814477664$ & $0.901920266712351$ \\ 
$0.3000$ & $1.339243905573495$ & $1.569217796713238$ \\ 
$0.5000$ & $1.057355422358385$ & $1.719659382215715$ \\ 
$0.6000$ & $0.991379731889011$ & $1.750119945541900$ \\ 
$0.8000$ & $1.119711839328917$ & $1.782048657856496$ \\ 
$1.0000$ & $1.404221050193842$ & $1.797421134420320$ \\
\hline
\end{tabular}
\end{center}
\caption{$L^{2}$-norm errors in MMAE and SCEM approximations for Example 1}
\label{table:2}
\end{table}

\begin{figure}
    \centering
    \includegraphics[scale=0.4]{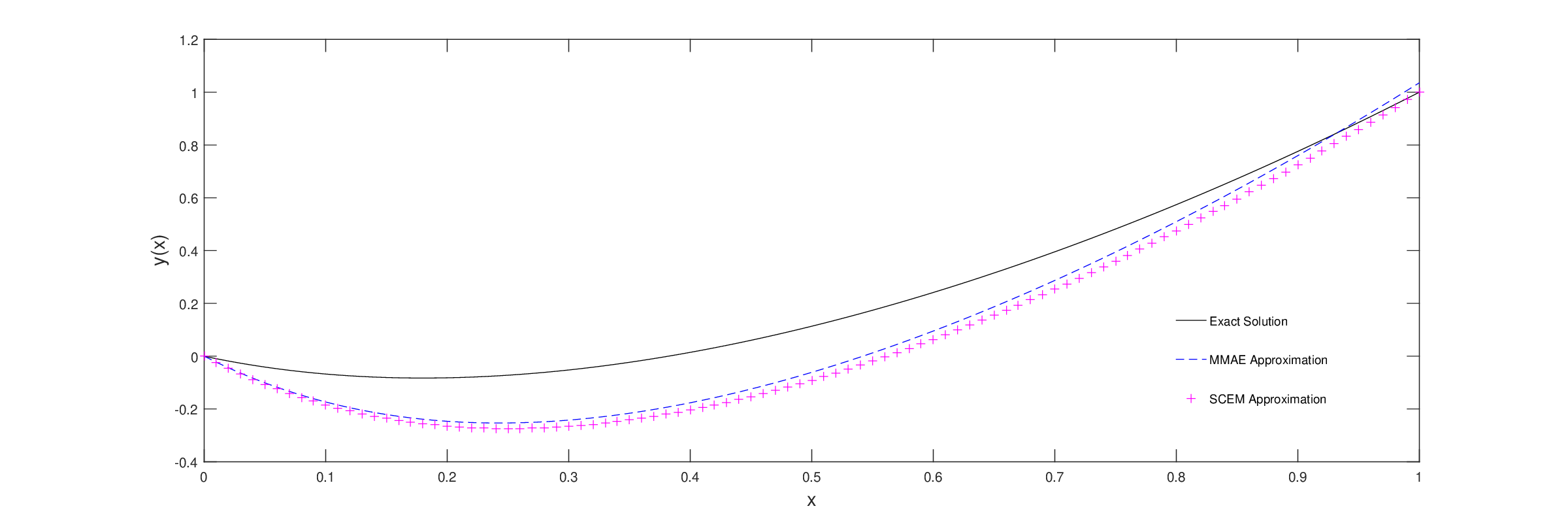}
    \caption{Comparison of results for Example 1, $\protect\varepsilon =0.3$}
    \label{fig:my_label11}
\end{figure}

\begin{figure}
    \centering
    \includegraphics[scale=0.4]{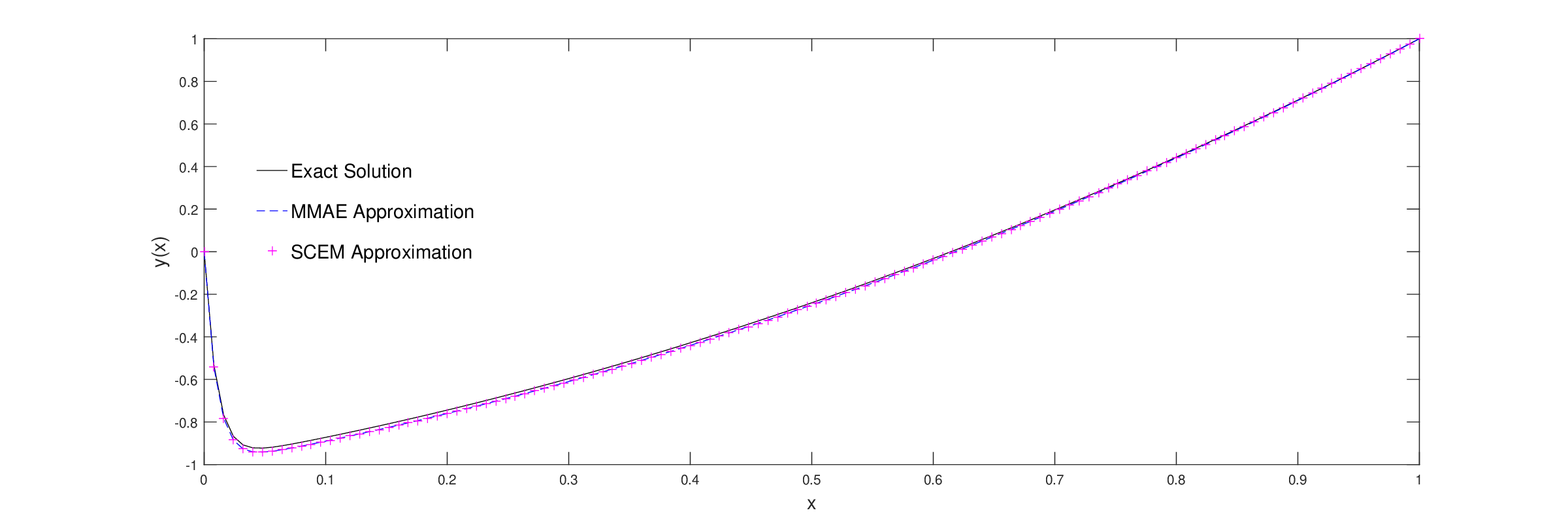}
    \caption{Comparison of results for Example 1, $\protect\varepsilon =0.01$}
    \label{fig:my_label12}
\end{figure}
Table \ref{table:2} shows that MMAE approximations are slightly more accurate than SCEM approximations for non-homogeneous problem, Example 1. As one can point out from Fig.\ref{fig:my_label11} and Fig.\ref{fig:my_label12}, as $\varepsilon \rightarrow 0^{+},$ MMAE
and SCEM approximations are getting more accurate. 

\textbf{Example 2:} Consider the singularly perturbed problem given in \cite{26}
\begin{equation}
-\varepsilon y^{\shortmid \shortmid }+y^{\shortmid }+y=1,\text{ }0\leq x<1,%
\text{ }y(0)=0,\text{ }y(1)=0.  \tag{3.15}
\end{equation}%
This problem has rapid changes near the point $x=1$. The exact solution to
(3.15) is given as

\begin{equation}
y(x)=1+\frac{e^{\lambda _{1}x}\left( e^{\lambda _{2}}-1\right) }{e^{\lambda
_{1}}-e^{\lambda _{2}}}-\frac{e^{\lambda _{2}x}\left( e^{\lambda
_{1}}-1\right) }{e^{\lambda _{1}}-e^{\lambda _{2}}},  \tag{3.16}
\end{equation}%
where $\lambda _{1}=\frac{1+\sqrt{1+4\varepsilon }}{2\varepsilon },$ $%
\lambda _{2}=\frac{1-\sqrt{1+4\varepsilon }}{2\varepsilon }$ and one-term
SCEM and MMAE approximations are obtained as follows.
The first term of the outer approximation is obtained letting $\varepsilon=0$ in the original equation (3.15)
\begin{equation}
y^{\shortmid }+y=1, \text{ }y_{0}(0)=0. \tag{3.17}
\end{equation}
It is obvious that the solution is 
\begin{equation}
    y_{0}(x)=1-e^{-x}. \tag{3.18}
\end{equation}
Since the problem (3.15) has boundary layer near the right end-point $x=1$,
the stretching variable should have a form  $\overline{x}=\frac{x-1}{%
\varepsilon }$.
The first term of the inner approximation, letting $\varepsilon =0$, is obtained as the solution of 
\begin{equation}
Y_{0}^{^{\shortparallel }}-Y_{0}^{^{\shortmid }}=0,\ Y_{0}(\overline{x}=0)=0\ \tag{3.19}
\end{equation}
and it is clear that the solution is
\begin{equation}
Y_{0}(\overline{x})=c_{1}(e^{\overline{x}}-1), \ c_{1}\in \mathbb{R}. \tag{3.20}
\end{equation}
Now that the first terms of outer and inner approximations are obtained, these approximations should be matched by the definition of MMAE.
\begin{equation}
\underset{x\longrightarrow 1^{-}}{\lim }{}y_{0}{(x)}=\underset{x\longrightarrow 1^{-}}{\lim }(1-e^{-x})=
\underset{\overline{x}\rightarrow -\infty }{%
\lim }c_{1}(e^{\overline{x}}-1)=\underset{\overline{x}\longrightarrow -\infty }{\lim }
Y_{0}(\overline{x}) \tag{3.21}
\end{equation}
Thus, it is determined that
\begin{equation}
Y_{0}(\overline{x})=\big( e^{-1}-1 \big)\big(e^{\overline{x}}-1\big). \tag{3.22}
\end{equation}
Finally, the composite MMAE approximation to the problem (3.15) is found as
\begin{equation}
y_{composite}(x,\varepsilon )=e^{-1}-e^{-x}+\left( e^{-1}-1\right) \left( e^{%
\frac{x-1}{\varepsilon }}-1\right). \tag{3.23}
\end{equation}
Now, the one-term SCEM approximation can be determined substituting the approximation 
\begin{equation}
y_{0}^{scem}(x,\overline{x},\varepsilon )=y_{0}(x,\varepsilon )+\Psi _{0}(%
\overline{x},\varepsilon )=1-e^{-x}+\Psi _{0}(%
\overline{x},\varepsilon ) \tag{3.24}
\end{equation}
into the original problem (3.15): 
\begin{align}
-\varepsilon \frac{d^2}{dx^2} \Big (1-e^{-x}+\Psi _{0}(%
\overline{x},\varepsilon )\Big )+ \frac{d}{dx}\Big(1-e^{-x}+\Psi _{0}(%
\overline{x},\varepsilon )\Big)+ \Big(1-e^{-x}+\Psi _{0}(%
\overline{x},\varepsilon) \Big)& =1 \notag \\  \notag
-\varepsilon \Big (-e^{-x}+\frac{\Psi_{0}^{\shortparallel}(\overline{x},\varepsilon )}{\varepsilon^{2}}\Big )+ \Big(e^{-x}+\frac{\Psi _{0} ^{\shortmid}(\overline{x},\varepsilon )}{\varepsilon}\Big)+\Big(1-e^{-x}+\Psi _{0}(%
\overline{x},\varepsilon) \Big)& =1\\ \tag{3.25}
\varepsilon ^{2} e^{-x} -\Psi _{0}^{\shortparallel}(\overline{x},\varepsilon )+\Psi _{0}^{\shortmid}(\overline{x},\varepsilon ) + \varepsilon ^{2}\Psi _{0} (\overline{x},\varepsilon ) & = 0.
\end{align}
In order to obtain first term of SCEM approximation, letting $\varepsilon=0$, the first complementary SCEM approximation $\Psi _{0}(%
\overline{x},\varepsilon )$ is obtained as the solution of the following two-point boundary value problem 
\begin{equation}
\Psi _{0}^{\shortparallel}(\overline{x},\varepsilon )-\Psi _{0}^{\shortmid}(\overline{x},\varepsilon )=0,\ \Psi_{0}\Big (0,\frac{-1}{\varepsilon },\varepsilon \Big ) =0,\ \Psi _{0} (1,0,\varepsilon )=e^{-1}-1. \tag{3.26}
\end{equation}
Finally, the first SCEM and MMAE approximation are found as follows

\begin{equation}
y_{0}^{scem}(x,\overline{x},\varepsilon )=1-e^{-x}+\left( e^{\frac{x-1}{%
\varepsilon }}-e^{\frac{-1}{\varepsilon }}\right) \left( \frac{1-\frac{1}{e}%
}{e^\frac{-1}{\varepsilon }-1}\right) ,  \tag{3.27}
\end{equation}

\begin{center}
\begin{equation}
y_{composite}(x,\varepsilon )=e^{-1}-e^{-x}+\left( e^{-1}-1\right) \left( e^{%
\frac{x-1}{\varepsilon }}-1\right) .  \tag{3.28}
\end{equation}
\end{center}

\begin{table}
\begin{center}
\begin{tabular}{|c|c|c|}
\hline
$\varepsilon $ & $L^{2}$ error in MMAE & $L^{2}$ error in SCEM \\ 
\hline
$0.0050$ & $0.138999385861808$ & $0.138999385861808$ \\ 
$0.0070$ & $0.192846875716269$ & $0.192846875716269$ \\ 
$0.0100$ & $0.271762063576098$ & $0.271762063576098$ \\ 
$0.0500$ & $1.121307087312208$ & $1.121307202107427$ \\ 
$0.0700$ & $1.413966007685341$ & $1.414000101581222$ \\ 
$0.1000$ & $1.703813183303206$ & $1.706213659924062$ \\ 
$0.2500$ & $1.129307299656383$ & $1.882651704001525$ \\ 
$0.5000$ & $4.188478595061564$ & $2.177262694969191$ \\ 
$0.7000$ & $7.946036759774072$ & $2.632927334749375$ \\ 
$0.8000$ & $9.570235541695139$ & $2.846650013414474$ \\ 
$0.9000$ & $11.023403578129694$ & $3.041045798307327$ \\ 
$1.0000$ & $12.321456684111308$ & $3.215647240288718$ \\
\hline
\end{tabular}
\end{center}
\caption{$L^{2}$-norm errors in MMAE and SCEM approximations for Example 2}
\label{table:3}
\end{table}

\begin{figure}
    \centering
    \includegraphics[scale=0.40]{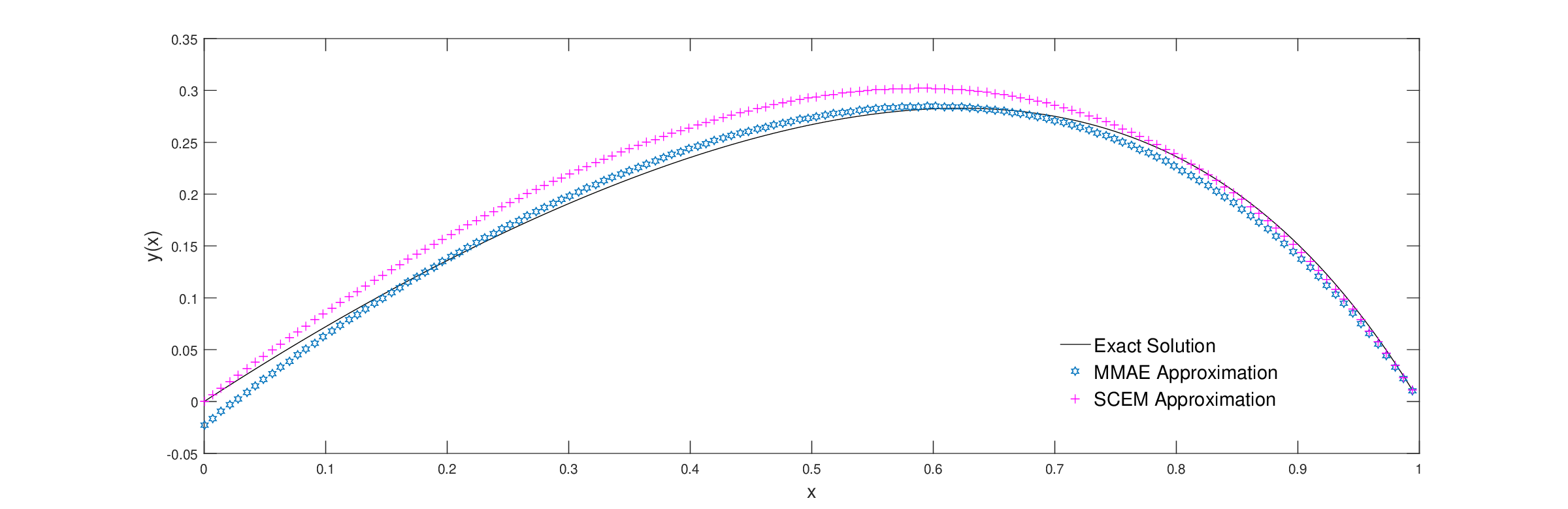}
    \caption{Comparison of results for Example 2, $\protect\varepsilon =0.3$}
    \label{fig:my_label21}
\end{figure}

\begin{figure}
    \centering
    \includegraphics[scale=0.4]{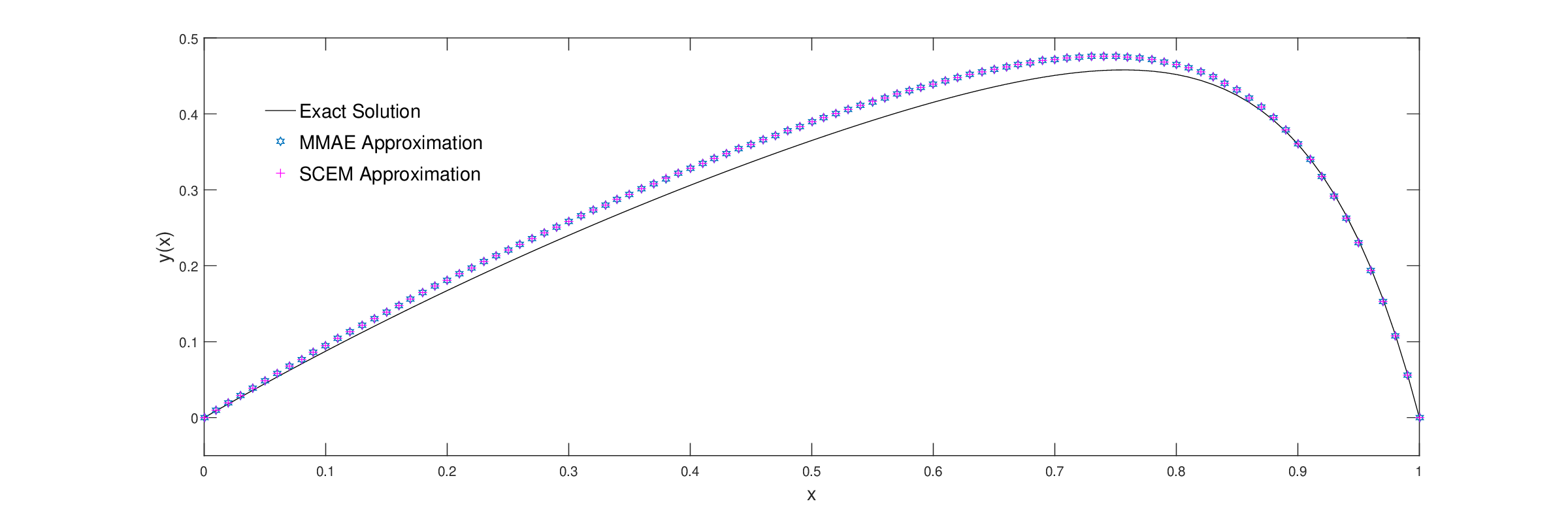}
    \caption{Comparison of results for Example 2, $\protect\varepsilon =0.1$}
    \label{fig:my_label22}
\end{figure}

\begin{figure}
    \centering
    \includegraphics[scale=0.4]{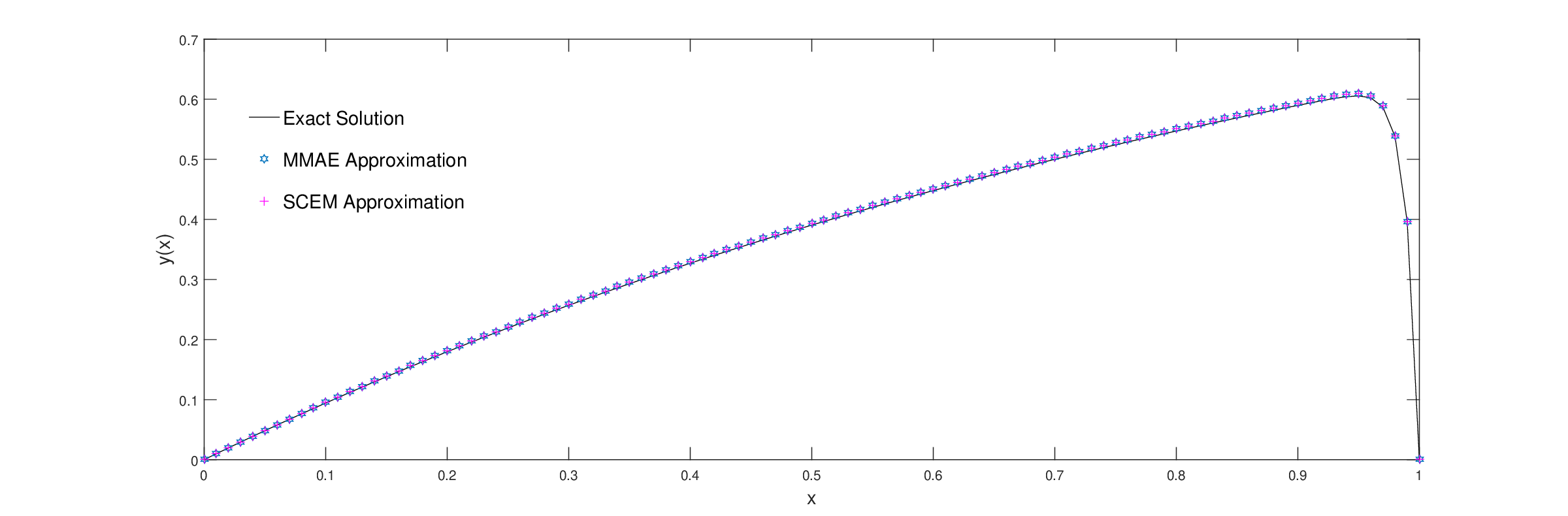}
    \caption{Comparison of results for Example 2, $\protect \varepsilon =0.01$}
    \label{fig:my_label23}
\end{figure}

In Fig.\ref{fig:my_label21} - Fig.\ref{fig:my_label23}, comparisons of approximations that are generated by SCEM and MMAE are given for $\varepsilon=0.3$, $\varepsilon=0.1$ and $\varepsilon=0.01$, respectively. It can be pointed out that for $\varepsilon=0.3$, MMAE generates relatively more accurate approximation. 
As it can be observed in Table \ref{table:3}, for more smaller $\varepsilon$ values, SCEM and MMAE approximations exactly overlap. In general manner, MMAE approximations are a bit more accurate than SCEM approximations, but for increasing $\varepsilon$ values, especially for larger values than $\varepsilon=0.5$, MMAE approximations are unacceptable. 

\section{Conclusion}

In this paper, the well-known method MMAE and relatively new one SCEM are
compared for solving singularly perturbed linear problems. An illustrative
example is given to demonstrate all the steps of both methods in detail. Two
non-homogeneous equations that has a left-end boundary layer and that has
right-end boundary layer are provided respectively in order to analyze
different cases. It is observed that SCEM gives more accurate approximations 
than MMAE approximations for homogeneous problems since non-homogeneous parts 
lead loss in captured terms during the balancing process. On the other hand, it is 
observed that MMAE approximations are a bit more accurate than those that are 
obtained by SCEM for solving non-homogeneous problems. 
Although only one-term approximations are proposed, SCEM and MMAE
gives highly accurate approximations for moderately small $\varepsilon $ values.
It is observed that if $\varepsilon $ is not kept sufficiently small, both SCEM and MMAE
approximations differ from the exact solutions.
The comparisons show that SCEM is an effective and flexible alternative method for
solving singularly perturbed, especially for homogeneous problems. Furthermore, 
the boundary conditions are satisfied exactly and any matching process is not required contrary to MMAE.

\section{Acknowledgments}

The author thanks anonymous reviewers for extremely helpful comments and suggestions to improve the quality of present paper. The author also thanks Professor Natesan Srinivasan of Indian Institute of Technology Guwahati for his helpful suggestions and recommendations about the paper.

The author dedicates his work to his grandfamilies S\"{u}leyman-Dudu Y\i lmaz and \c{S}emsi-Ay\c{s}e Cengizci.

\end{document}